\documentclass[reqno]{amsart} \usepackage{amscd}
 \usepackage{epsf}
 \usepackage{epsfig}
 \usepackage{graphics}
 \usepackage{graphicx}
 \usepackage{tikz} 
 \usepackage{amssymb}
 \usepackage{tkz-euclide}
 \usepackage{verbatim}
 \usepackage[english]{babel}
 \usepackage{tabu}
 \usepackage{amsmath}

 \theoremstyle{remark} 
 \numberwithin{equation}{section}

 \providecommand \hf{\hspace*{0.5cm}}
 
 \theoremstyle{plain}
 \newtheorem{thm}{Theorem}[section]
 
 \newtheorem{lmm}[thm]{Lemma}

 \newtheorem{rem}[thm]{Remark}
 \newtheorem{mresult}[thm]{Main Result}

 \theoremstyle{definition}

 \usepackage{hyperref}

\begin{document}
 	\title{Enumeration of rational curves in a moving family of $\mathbb{P}^2$}
 	%

 	\author[R. Mukherjee]{Ritwik Mukherjee}
 	
 	\address{School of Mathematics, National Institute of Science Education and Research, Bhubaneswar (HBNI), Odisha 752050, India}
 	
 	\email{ritwikm@niser.ac.in}
 	
 	\author[Anantadulal Paul]{Anantadulal Paul}
 	
 	\address{School of Mathematics, National Institute of Science Education and Research, Bhubaneswar (HBNI), Odisha 752050, India}
 	
 	\email{anantadulal.paul@niser.ac.in } 
 	
 	\author[Rahul Kumar Singh]{Rahul Kumar Singh}
 	
 	\address{School of Mathematics, National Institute of Science Education and Research, Bhubaneswar (HBNI), Odisha 752050, India}
 	
 	\email{rahul412@niser.ac.in }


 	\begin{abstract}
 		We obtain a recursive formula for the number of rational degree $d$ curves in $\mathbb{P}^3$, whose image lies in a $\mathbb{P}^2$, passing through $r$ 
 		lines and $s$ points, where $r + 2s = 3d+2$. This can be viewed as a family version of the classical question of counting rational curves in $\mathbb{P}^2$. 
 		We verify that our numbers are consistent with those obtained by T.~Laarakker, where he studies the parallel question of counting $\delta$-nodal degree $d$ 
 		curves in $\mathbb{P}^3$ whose image lies inside a $\mathbb{P}^2$. Our numbers give evidence to support the conjecture, that the polynomials obtained by 
 		T.~Laarakker 
 		are enumerative when $d \geq 1 + [\frac{\delta}{2}]$, which is analogous to the {G}\"ottsche threshold for counting nodal curves in $\mathbb{P}^2$. 
 	\end{abstract}
 	
 	\maketitle
 	
 	\tableofcontents

 	\section{Introduction}
 	One of the most fundamental and studied problems in enumerative geometry is the following: what is $N_d^{\delta}$, the number of 
 	degree $d$ curves in $\mathbb{P}^2$ that have $\delta$ distinct nodes and pass through $\frac{d(d+3)}{2} - \delta$ generic points?
 	The question was studied more than a hundred years ago by Zeuthen (\cite{Zeuthen}) and has been studied extensively in the last thirty 
 	years 
 	by Ran (\cite{Ran1}, \cite{Ran2}), Vainsencher (\cite{Van}), Caporaso-Harris (\cite{CH}), Kazarian (\cite{Kaz}), Kleiman and Piene (\cite{KP1}), 
 	Florian Block (\cite{FB}),
 	Tzeng and Li (\cite{Tz}, \cite{Tzeng_Li}), Kool, Shende and Thomas (\cite{KST}), Berczi (\cite{Berczi}) and Fomin and Mikhalkin (\cite{FoMi}), 
 	amongst others. This question has been investigated from several perspectives and is very well understood. \\ 
 	\hf \hf The problem motivates a natural generalization considered by Kleiman and Piene in \cite{KP2}, where they study the 
 	enumerative geometry of singular curves in a moving family of 
 	surfaces. More recently, this question has been studied further by  T.~Laarakker in \cite{TL}, where he obtains 
 	a formula for the following number: how many degree $d$ curves are 
 	there in $\mathbb{P}^3$ whose image lies in a $\mathbb{P}^2$, that pass through $\frac{d(d+3)}{2} + 3 - \delta$ generic lines and have $\delta$-nodes 
 	(provided $d\geq \delta$). This can be viewed as a family version of the classical problem of computing $N_d^{\delta}$. \\ 
 	\hf \hf Motivated by the papers of Kleiman and Piene (\cite{KP2}) and T.~Laarakker (\cite{TL}), 
 	we have studied the parallel question of counting stable 
 	rational maps into a family of moving target spaces. This can be viewed as a family version of the famous question of enumerating rational curves 
 	in $\mathbb{P}^2$, that was studied by Kontsevich-Manin (\cite{KoMa}) and Ruan-Tian (\cite{RT}). The main result of our paper is as follows: 
 	
 	\begin{mresult}
 		Let $N_{d}^{\mathbb{P}^3, \mathrm{Planar}}(r,s)$ be the number of genus zero, degree $d$ curves in $\mathbb{P}^3$, whose image  lies in a $\mathbb{P}^2$, intersecting 
 		$r$ generic lines and $s$ generic points (where $r  + 2s = 3d+2$). We have 
 		a recursive formula to compute $N_{d}^{\mathbb{P}^3, \mathrm{Planar}}(r,s)$ for all $d \geq 2$.
 	\end{mresult}
 	
 	\begin{rem}
 		Note that for $d=1$, the corresponding question is classical Schubert calculus and there $r+2s = 4$ as opposed to $5$. 
 	\end{rem}
 	
 	\begin{rem}
 		We note that when $s=3$ and $d \geq 2$, the number $N_{d}^{\mathbb{P}^3, \mathrm{Planar}}(3d-4,3)$ is the number of rational curves in $\mathbb{P}^2$ 
 		through $3d-1$ points; this is because $3$ generic  points in $\mathbb{P}^3$ determine a unique $\mathbb{P}^2$. We also note that when $s>3$, 
 		$N_d^{\mathbb{P}^3, \mathrm{Planar}}(r,s)$ is zero, since $4$ or more generic points do not lie in a plane.  
 	\end{rem}

 	We have written a mathematica program to implement our formula; 
 	the program is available on our web page
 	\[ \textnormal{\url{https://www.sites.google.com/site/ritwik371/home}}. \]
 	In section \ref{low_deg_check}, we verify that the numbers we compute are logically consistent with those obtained by T.~Laarakker in \cite{TL} till $d=6$. 
 	This gives strong evidence to support the conjecture that his formulas for $\delta$-nodal planar degree $d$ curves in $\mathbb{P}^3$ 
 	are expected to  be enumerative when $d \geq 1 + [\frac{\delta}{2}]$ 
 	(as opposed to $d \geq \delta$ which is proved in \cite{TL}). 
 	Starting from $d=7$, we can not use the result of \cite{TL} to make any consistency check, since the corresponding nodal polynomial is not expected to 
 	be enumerative (due to the presence of double lines); this is explained in section \ref{low_deg_check}.  
 	

 	
 	
 	
 	
 	
 	\section{Notation }
 	Let us define a \textbf{planar} curve in $\mathbb{P}^3$ to be a curve, whose image lies inside a $\mathbb{P}^2$. 
 	We will now develop some notation to describe the space of planar curves of a given degree $d$. \\ 
 	\hf \hf Let us denote the dual of $\mathbb{P}^3$ by 
 	$\hat{\mathbb{P}}^3$; this is the space of $\mathbb{P}^2$ inside $\mathbb{P}^3$. 
 	An element of $\hat{\mathbb{P}}^3$ can be thought of as 
 	a non zero linear functional $\eta: \mathbb{C}^4 \longrightarrow \mathbb{C}$ upto scaling (i.e., it is the projectivization of the dual of $\mathbb{C}^4$).  
 	Given such an $\eta$, we  define the projectivization of its zero set as $\mathbb{P}^2_{\eta}$. In other words, 
 	\begin{align*}
 	\mathbb{P}^2_{\eta} &:= \mathbb{P}(\eta^{-1}(0)). 
 	\end{align*}
 	Note that this $\mathbb{P}^2_{\eta}$ is a subset of $\mathbb{P}^3$. Next, when $d \geq 2$, we  
 	define the moduli space of planar degree $d$ curves into $\mathbb{P}^3$ 
 	as a fibre bundle over $\hat{\mathbb{P}}^3$. More precisely, 
 	we define 
 	\[ \pi: \overline{\mathcal{M}}_{0,k}^{\mathrm{Planar}}(\mathbb{P}^3,d) \longrightarrow \hat{\mathbb{P}}^3 \] 
 	to be the fiber bundle, such that 
 	\[ \pi^{-1}([\eta]) := \overline{\mathcal{M}}_{0,k}(\mathbb{P}^2_{\eta}, d). \]
 	Here we are using the standard notation to denote $\mathcal{M}_{0,k}(X, \beta)$ to be the moduli space of genus zero stable maps, 
 	representing the class $\beta \in H_2(X, \mathbb{Z})$ and $\overline{\mathcal{M}}_{0,k}(X, \beta)$ to be its stable map compactification. 
 	Since the dimension of a fiber bundle is dimension of the base, plus the dimension of the fiber, we conclude that the dimension of 
 	$\overline{\mathcal{M}}_{0,k}^{\mathrm{Planar}}(\mathbb{P}^3,d)$ is $3d+2+k$.\\
 	\hf\hf Next, we note that the space of planes in $\mathbb{P}^3$ can also be thought of as the Grassmanian $\mathbb{G}(3,4)$.  Let 
 	$\gamma_{3,4}$ denote the tautological three plane bundle over the Grassmanian. Since $\mathbb{G}(3,4)$ can be identified with $\hat{\mathbb{P}}^3$, 
 	we can think of $\gamma_{3,4}$ as a bundle over $\hat{\mathbb{P}}^3$. \\ 
 	\hf \hf When $d=1$, we define $\overline{\mathcal{M}}_{0,0}^{\mathrm{Planar}}(\mathbb{P}^3,1)$ to be 
 	\[ \overline{\mathcal{M}}_{0,0}^{\mathrm{Planar}}(\mathbb{P}^3,1) := \mathbb{P}(\gamma_{3,4}^*) \longrightarrow \hat{\mathbb{P}}^3.\] 
 	We note that an element of $\overline{\mathcal{M}}_{0,0}^{\mathrm{Planar}}(\mathbb{P}^3,1)$ 
 	is of the form $(L, H)$, where $L$ is a line in $\mathbb{P}^3$ and $H$ is a plane containing $L$. 
 	Since a line is not  contained in a unique plane, $\overline{\mathcal{M}}_{0,0}^{\mathrm{Planar}}(\mathbb{P}^3,1)$ is not the 
 	same as the space of lines; infact we note that the space of  lines is $4$ 
 	dimensional, while the dimension of $\overline{\mathcal{M}}_{0,0}^{\mathrm{Planar}}(\mathbb{P}^3,1)$ is $5$.\\  
 	\hf \hf We will now define a few numbers by intersecting cycles on $\overline{\mathcal{M}}_{0,0}^{\mathrm{Planar}}(\mathbb{P}^3,d)$, the moduli space with zero marked  points 
 	(this includes the case $d=1$; unless otherwise stated we always include the case $d=1$ in any of our statements). 
 	Let 
 	$\mathcal{H}_L$ and $\mathcal{H}_p$  denote the classes  of the cycles in $\overline{\mathcal{M}}_{0,0}^{\mathrm{Planar}}(\mathbb{P}^3,d)$ that corresponds 
 	to the subspace of curves passing through a generic line and a point respectively. 
 	We also denote  $a$ 
 	to be the standard generator of $H^*(\hat{\mathbb{P}}^3; \mathbb{Z})$.  
 	Let us now define 
 	\begin{align}
 	N_d^{\mathbb{P}^3,\mathrm{Planar}}(r,s,\theta) &:= \langle a^{\theta}, ~~\overline{\mathcal{M}}_{0,0}^{\mathrm{Planar}}(\mathbb{P}^3,d) \cap \mathcal{H}_L^r \cap \mathcal{H}_p^s \rangle.  
 	\label{nd_rs_theta}
 	\end{align}
 	We formally define the the left hand side of $\eqref{nd_rs_theta}$ to be zero if $r+2s + \theta \neq 3d+2$ 
 	(since the right hand side of \eqref{nd_rs_theta} doesn't make sense unless $r+2s + \theta = 3d+2$).     
 	Note that when $\theta =0$, $r + 2s =3d+2$ and $d \geq 2$,  
 	$N_d^{\mathbb{P}^3,\mathrm{Planar}}(r,s,0)$ is precisely equal to the number of rational planar degree $d$-curves in $\mathbb{P}^3$ intersecting $r$ generic lines and 
 	$s$ generic points. \\ 
 	\hf \hf Next, we will define a number $B_{d_1,d_2}(r_1,s_1,r_2,s_2,\theta)$ by intersecting it on the product 
 	of two  moduli spaces as follows: 
 	\begin{eqnarray}
 	B_{d_1,d_2}(r_1,s_1,r_2,s_2,\theta) &:= \Big\langle \pi_1^*( \mathcal{H}_L^{r_1} \mathcal{H}_p^{s_1})\cdot \pi_2^*( \mathcal{H}_L^{r_2} \mathcal{H}_p^{s_2}) (\pi_2^*a)^{\theta}\cdot 
 	(\pi^*\Delta), \nonumber \\ 
 	& \qquad [\overline{\mathcal{M}}_{0,0}^{\mathrm{Planar}}(\mathbb{P}^3,d_1) \times \overline{\mathcal{M}}_{0,0}^{\mathrm{Planar}}(\mathbb{P}^3,d_2)] \Big\rangle.  \label{bd_delta}
 	\end{eqnarray}
 	Here $\Delta$ denotes the class of the diagonal in $\hat{\mathbb{P}}^3 \times \hat{\mathbb{P}}^3$ and $\pi_1$ and $\pi_2$ are the obvious projection 
 	maps. 
 	Again, we formally define the left hand 
 	side of \eqref{bd_delta} to be zero, unless $r_1 +  2s_1 + r_2 + 2 s_2 + \theta = 3d_1 + 3d_2 + 4$ (since in that case, 
 	the right hand side doesn't  make sense). Note that when $\theta =0$ and $r_1 +  2s_1 + r_2 + 2 s_2 = 3d_1 + 3d_2 + 4$, the number $B_{d_1,d_2}(r_1,s_1,r_2,s_2,0)$ 
 	denotes the number of two component rational curves in $\mathbb{P}^3$ of degree $(d_1, d_2)$, such that the two components intersect and lie in a common plane  
 	and the $d_i$ component intersects $r_i$ lines and $s_i$ points, for $i=1$ and $2$.  
 	
 	\section{Recursive Formula and its Proof}
 	We are now ready to state our recursion formula 
 	\begin{lmm}
 		\label{base_case_recursion}
 		If $d=1 $, then the number $N_d^{\mathbb{P}^3, \mathrm{Planar}}(r,s,\theta)$ is given by 
 		\begin{align}
 		\label{base_formula}
 		N_1^{\mathbb{P}^3,\mathrm{Planar}}(r,s,\theta) & = \begin{cases} 0 & \mbox{if} ~~(r,s,\theta) =(1,2,0), \\ 
 		0 & \mbox{if }  ~~(r,s,\theta) =(3,1,0),\\ 
 		0 & \mbox{if } ~~(r,s,\theta) =(5,0,0),\\ 
 		1 &  \mbox{if} ~~(r,s,\theta) =(0,2,1),\\ 
 		1 & \mbox{if} ~~(r,s,\theta) =(2,1,1),\\ 
 		2 & \mbox{if} ~~(r,s,\theta) =(4,0,1), \\ 
 		1 & \mbox{if} ~~(r,s,\theta) =(1,1,2), \\ 
 		2 & \mbox{if} ~~(r,s,\theta) =(3,0,2), \\ 
 		0 & \mbox{if} ~~(r,s,\theta) =(0,1,3), \\ 
 		1 & \mbox{if} ~~(r,s,\theta) =(2,0,3), \\
 		0 & \mbox{\textnormal{otherwise}}.\end{cases} 
 		\end{align}
 	\end{lmm}
 	
 	\begin{lmm}
 		\label{second_case_recursion}
 		If $d=2 $, then the number $N_d^{\mathbb{P}^3, \mathrm{Planar}}(r,s,\theta)$ is given by 
 		\begin{align}
 		\label{second_case_formula}
 		N_2^{\mathbb{P}^3,\mathrm{Planar}}(r,s,\theta) & = \begin{cases} 92 & \mbox{if} ~~(r,s,\theta) =(8,0,0), \\ 
 		18 & \mbox{if }  ~~(r,s,\theta) =(6,1,0),\\ 
 		4 & \mbox{if } ~~(r,s,\theta) =(4,2,0),\\ 
 		1 &  \mbox{if} ~~(r,s,\theta) =(2,3,0),\\ 
 		34 & \mbox{if} ~~(r,s,\theta) =(7,0,1),\\ 
 		6 & \mbox{if} ~~(r,s,\theta) =(5,1,1), \\ 
 		1 & \mbox{if} ~~(r,s,\theta) =(3,2,1), \\ 
 		0 & \mbox{if} ~~(r,s,\theta) =(1,3,1), \\ 
 		8 & \mbox{if} ~~(r,s,\theta) =(6,0,2), \\ 
 		1 & \mbox{if} ~~(r,s,\theta) =(4,1,2), \\ 
 		0 & \mbox{if} ~~(r,s,\theta) =(2,2,2), \\
 		0 & \mbox{if} ~~(r,s,\theta) =(0,3,2), \\
 		1 & \mbox{if} ~~(r,s,\theta) =(5,0,3), \\
 		0 & \mbox{if} ~~(r,s,\theta) =(3,1,3), \\
 		0 & \mbox{if} ~~(r,s,\theta) =(1,2,3), \\
 		0 & \mbox{\textnormal{otherwise}}.\end{cases} 
 		\end{align}
 	\end{lmm}

 	%
 	
 	
 	\begin{thm}
 		\label{rec_formula_mthm}
 		If $d\geq 3$, then 
 		\begin{align}
 		\label{triv_eqn}
 		N_d^{\mathbb{P}^3, \mathrm{Planar}}(r,s,\theta) &=  \begin{cases} 0 & \mbox{if} ~~r + 2s + \theta \neq 3d+2, \\ 
 		0 & \mbox{if }  ~~s>3,\\ 
 		0 & \mbox{if } ~~\theta>3.\end{cases} 
 		\end{align}
 		In the remaining case  when $r+2s + \theta = 3d+2$, $s \leq 3$ and $\theta \leq 3$, we have   
 		
 		\begin{align}\label{recurfor}
 		N_d^{\mathbb{P}^3, \mathrm{Planar}}(r,s,\theta)& =2d N_d^{\mathbb{P}^3, \mathrm{Planar}}(r-2,s+1,\theta) \nonumber \\ 
 		& +\sum_{r_1 =0}^{r-3}\sum_{s_1 =0}^{s}\sum_{d_1=1}^{d-1}\binom{r-3}{r_1}\binom{s}{s_1}d_1^2 d_2\times \nonumber \\ 
 		& \Big(d_2 B_{d_1,d_2}(r_1+1,s_1,r_2-2,s_2,\theta)
 		-d_1B_{d_1,d_2}(r_1,s_1,r_2-1,s_2,\theta)\Big), 
 		\end{align}
 		where in the above expression, $d_2:= d-d_1$, $r_2:= r-r_1$ and $s_2:= s-s_1$. 
 		Furthermore, 
 		$\forall d_1, d_2 \geq 1$,  we have  \\
 		\begin{equation}\label{b1}
 		B_{d_1,d_2}(r_1,s_1,r_2,s_2,\theta) =\sum_{i=0}^3 N_{d_1}^{\mathbb{P}^3, \mathrm{Planar}}(r_1,s_1,i)\times N_{d_2}^{\mathbb{P}^3, \mathrm{Planar}}(r_2,s_2,\theta+3-i).
 		\end{equation}
 	\end{thm}
 	
 	\begin{rem}
 		We note that equations \eqref{base_formula}, \eqref{second_case_formula}, \eqref{triv_eqn},  
 		\eqref{recurfor} and \eqref{b1} allow us to compute $N_d^{\mathbb{P}^3, \mathrm{Planar}}(r,s,\theta)$ 
 		for all $d \geq 1$ and all $r, s, \theta \geq 0$.
 	\end{rem}
 	
 	
 	\noindent \textbf{Proof of Theorem \ref{rec_formula_mthm}:} 
 	We will start by proving equation \eqref{triv_eqn}. The first equation is true simply because we are pairing a cohomology class with a homology 
 	class of different dimensions (see the remark after equation \eqref{nd_rs_theta}). Next, when $s>3$, we note that there can not be any planar curves, 
 	because $4$ generic points do not lie on a plane. Finally, we  note that $a^{\theta} =0 \in H^*(\hat{\mathbb{P}}^3;\mathbb{Z})$ if $\theta >3$, which proves the last equality. \\ 
 	\hf \hf We now justify the main thing, which is equation \eqref{recurfor}.
 	The idea is very similar to the idea used to compute 
 	the number of rational degree $d$ curves in $\mathbb{P}^2$ (and also $\mathbb{P}^3$) that is given in 
 	\cite{KoMa}, \cite{McSa} and \cite{FuPa}.
 	As in the case of counting curves in $\mathbb{P}^2$,  
 	let us first consider $\overline{\mathcal{M}}_{0,4}$. This space is isomorphic to $\mathbb{P}^1$; hence we have equivalence of the 
 	following divisors   
 	\begin{align}
 	(x_1x_2|x_3 x_4) \equiv (x_1x_3|x_2 x_4). \label{wdvv_p1}
 	\end{align}
 	Here $(x_1x_2|x_3 x_4)$ denotes the domain which is a wedge of two spheres with the marked points $x_1$ and $x_2$ on the first 
 	sphere and $x_3$ and $x_4$ on the second sphere. The domain $(x_1x_3|x_2 x_4)$ is defined similarly. Since $\overline{\mathcal{M}}_{0,4} \approx \mathbb{P}^1$ 
 	is path connected, any two points determine the same divisor. \\
 	\hf \hf Let us now consider the projection map 
 	\[\pi: \overline{\mathcal{M}}_{0,4}^{\textrm{Planar}}(\mathbb{P}^3, d) \longrightarrow  \overline{\mathcal{M}}_{0,4}. \] 
 	We define a cycle $\mathcal{Z}$ in $\overline{\mathcal{M}}_{0,4}^{\textrm{Planar}}(\mathbb{P}^3, d)$, given by 
 	\begin{align*}
 	\mathcal{Z} &:= \textnormal{ev}_1^*([H]) \cdot \textnormal{ev}_2^*([H]) \cdot \textnormal{ev}_3^*([L]) \cdot \textnormal{ev}_4^*([L]) \cdot \mathcal{H}_L^{r-3} 
 	\cdot \mathcal{H}_p^s \cdot a^{\theta},  
 	\end{align*}
 	where $[H]$ and $[L]$ denote the class of a hyperplane and a line in $\mathbb{P}^3$, $\textnormal{ev}_i$ denotes the evaluation map at the $i^{\textnormal{th}}$ 
 	marked point and $a$ denotes the generator of $H^*(\hat{\mathbb{P}}^3;\mathbb{Z})$. 
 	
 	Let us now intersect the cycle $\mathcal{Z}$ by pulling back the left hand side of \eqref{wdvv_p1}, via $\pi$. 
 	We now note that 
 	\begin{align}
 	\pi^*(x_1 x_2|x_3 x_3) \cdot \mathcal{Z}  &=  
 	N_d^{\mathbb{P}^3,\mathrm{Planar}}(r,s,\theta)\; \nonumber \\ 
 	& +\sum_{r_1=0}^{r-3}\sum_{s_1=0}^s\sum_{d_1=1}^{d-1}\binom{r-3}{r_1}\binom{s}{s_1}d_1^3 d_2\times B_{d_1,d_2}(r_1,s_1,r_2-1,s_2,\theta),\label{wdvv1}
 	\end{align}
 	where in the above expression, $d_2:= d-d_1$, $r_2:= r-r_1$ and $s_2:= s-s_1$. \\ 
 	\hf \hf Next, we note that 
 	\begin{align}
 	\pi^*(x_1 x_3|x_2 x_4) \cdot \mathcal{Z}  &=  2d N_d^{\mathbb{P}^3,\mathrm{Planar}}(r-2,s+1,\theta)\; \nonumber \\ 
 	& +\sum_{r_1=0}^{r-3}\sum_{s_1=0}^s\sum_{d_1=1}^{d-1}\binom{r-3}{r_1}\binom{s}{s_1}d_1^2 d_2^2 \times B_{d_1,d_2}(r_1+1,s_1,r_2-2,s_2,\theta),
 	\label{wdvv2}
 	\end{align}
 	where as before $d_2:= d-d_1$, $r_2:= r-r_1$ and $s_2:= s-s_1$. \\ 
 	\hf \hf We now note that equations \eqref{wdvv_p1}, \eqref{wdvv1} and \eqref{wdvv2}, imply our desired recursive formula \eqref{recurfor}. \\  
 	\hf \hf Next, we will justify the formula for $B_{d_1, d_2}(r_1,  s_1, r_2, s_2,\theta)$ (equation \eqref{b1}). 
 	This follows immediately from the fact that the class of the diagonal 
 	is given by 
 	\begin{align*}
 	\Delta_{\hat{\mathbb{P}}^3 \times \hat{\mathbb{P}}^3} &= \pi_1^* a^3 + \pi_1^*a^2 \cdot \pi_2^*a + \pi_1^*a \cdot \pi_2^*a^2 +  \pi_2^* a^3,  
 	\end{align*}
 	where $a$ denotes the generator of $H^*(\hat{\mathbb{P}}^3;\mathbb{Z})$ and $\pi_1, \pi_2$ denote the respective projection maps. The formula now follows immediately 
 	from the definition of $$B_{d_1, d_2}(r_1,  s_1, r_2, s_2,\theta)$$  (namely, equation \eqref{bd_delta}).\\ 
 	\hf \hf It remains to prove the two base cases of the recursion, namely Lemma \ref{base_formula} and Lemma \ref{second_case_formula}.  \qed \\

 	\noindent \textbf{Proof of Lemma \ref{base_case_recursion}:} We recall that $\overline{\mathcal{M}}_{0,0}^{ \mathrm{Planar}}(\mathbb{P}^3, 1)$ 
 	is defined to be the projective bundle $\mathbb{P}(\gamma_{3,4}^*) \longrightarrow \mathbb{\hat{\mathbb{P}}}^3$. Now, we note that the 
 	Chern classes of the rank three vector bundle $\gamma_{3,4}^* \longrightarrow \mathbb{\hat{\mathbb{P}}}^3$ are given by 
 	\begin{align*}
 	c_i(\gamma_{3,4}^*) &= a^i \in H^{2i} (\mathbb{\hat{\mathbb{P}}}^3; \mathbb{Z}).  
 	\end{align*} 
 	The reason is explained in \cite[Page 18]{Zing_Notes}. Here $a$ is the standard generator of 
 	$H^{*} (\mathbb{\hat{\mathbb{P}}}^3; \mathbb{Z})$.\\ 
 	\hf \hf Next, we note that $a^i =0$ if $i >3$. Hence, the cohomology ring of $H^*(\mathbb{P}(\gamma_{3,4}^*))$ is given by 
 	\begin{align}
 	H^*(\mathbb{P}(\gamma_{3,4}^*)) &\approx \frac{\mathbb{Z}[a, \lambda]}{\langle \lambda^{3} + \lambda ^2 a + \lambda a^2 + a^3\rangle } 
 	\end{align}
 	where $\tilde{\gamma}\longrightarrow \mathbb{P}(\gamma_{3,4}^*)$ 
 	is the tautological line bundle over the projectivized bundle $\mathbb{P}(\gamma_{3,4}^*)$
 	and $\lambda := c_1(\tilde{\gamma}^*) \in H^*(\mathbb{P}(\gamma_{3,4}^*))$. This follows from \cite[Page 270]{BoTu}.\\ 
 	\hf \hf We now note the following two important facts: 
 	intersecting a generic line, corresponds to the cycle 
 	\[ \mathcal{H}_l = \lambda + a.\] 
 	Furthermore, passing through a generic point, corresponds to the cycle 
 	\[ \mathcal{H}_p = \lambda a. \] 
 	The reason for this can again be found in \cite[Pages 18 and 19]{Zing_Notes}. Hence, 
 	to compute $N_1^{\mathbb{P}^3, \mathrm{Planar}}(r, s, \theta)$ we have to compute the expression 
 	\[ (\lambda+a)^r (\lambda a)^s a^{\theta}, \] 
 	use the relationship 
 	\[ \lambda^3 = -(\lambda ^2 a + \lambda a^2 + a^3)  \]
 	and extract the coefficient of $\lambda^2 a^3$. This gives us all the numbers for various values of $r, s$ and $\theta$. \qed \\ 
 	
 	\noindent \textbf{Proof of Lemma \ref{second_case_recursion}:} First we note that every conic in $\mathbb{P}^3$ lies inside a unique plane 
 	(except a double line). 
 	Hence, let us consider the projective bundle
 	\[\mathbb{P}(\textnormal{Sym}^2(\gamma_{3,4}^*)) \longrightarrow \mathbb{\hat{\mathbb{P}}}^3.\] 
 	This space $\mathbb{P}(\textnormal{Sym}^2(\gamma_{3,4}^*))$ is the space of conics in $\mathbb{P}^3$ and a plane 
 	that contains the conic. The space of all smooth conics form an open dense subspace of $\mathbb{P}(\textnormal{Sym}^2(\gamma_{3,4}^*))$. 
 	Hence, to compute the numbers $N_2^{\mathbb{P}^3, \mathrm{Planar}}(r,s, \theta)$ (which is defined as an intersection number on 
 	$\overline{\mathcal{M}}^{\mathrm{Planar}}_{0,0}(\mathbb{P}^3,2)$), we can instead compute the relevant intersection number on 
 	$\mathbb{P}(\textnormal{Sym}^2(\gamma_{3,4}^*))$.\\ 
 	\hf \hf Next, we note that $\mathbb{P}(\textnormal{Sym}^2(\gamma_{3,4}^*))$ 
 	is a $\mathbb{P}^5$ bundle over $\mathbb{\hat{\mathbb{P}}}^3$. 
 	The cohomology ring structure of 
 	the total space is given by
 	\begin{align}
 	\label{ring_conic}
 	H^*(\mathbb{P}(\textnormal{Sym}^2(\gamma_{3,4}^*)) &\approx \frac{\mathbb{Z}[a, \lambda]}{\langle \lambda^6 + 4 \lambda^5 a + 10 \lambda^4 a^2 + 20\lambda^3 a^3 \rangle}. 
 	\end{align}
 	This follows from \textit{splitting principle} (see page $275$ in \cite{BoTu}).
 	We now note the following two important facts: 
 	intersecting a generic line, corresponds to the cycle 
 	\[ \mathcal{Z}_l = \lambda + 2a.\] 
 	Furthermore, passing through a generic point, corresponds to the cycle 
 	\[ \mathcal{Z}_p = \lambda a. \] 
 	The reason for this can again be found in \cite[Pages 18 and 19]{Zing_Notes}.
 	Hence, to compute $N_2^{\mathbb{P}^3, \mathrm{Planar}}(r, s, \theta)$ we have to compute the expression 
 	\[ (\lambda+2a)^r (\lambda a)^s a^{\theta}, \] 
 	use the relationship given by the cohomology ring structure in \eqref{ring_conic}, i.e.  
 	\[ \lambda^6 = -(4 \lambda^5 a + 10 \lambda^4 a^2 + 20\lambda^3 a^3)   \]
 	and extract the coefficient of $\lambda^5 a^3$. This gives us all the numbers for various values of $r, s$ and $\theta$. \qed \\

 	\section{Low degree checks} 
 	\label{low_deg_check}
 	We now describe concrete low degree checks that we have performed. Using our recursive formula, we have obtained the following numbers:
 	\begin{center}
 		\begin{tabular}{|c|c|c|c|c|} 
 			\hline
 			$(d,r,s)$ & $(3,11,0)$ & $(4,14,0)$ & $(5,17,0)$ & $(6,20,0)$ \\ 
 			\hline
 			$N_d^{\mathbb{P}^3, \mathrm{Planar}}(r,s)$ & $12960$ & $3727920$ & $1979329280$ & $1763519463360$ \\ 
 			\hline
 		\end{tabular}
 	\end{center}
 	Next, let $N^{\textnormal{Node}, \delta}_d(r,s)$ denote the number of planar degree $d$ curves in $\mathbb{P}^3$ with $\delta$ (unordered) nodes 
 	intersecting $r$ generic lines and $s$ generic points. These numbers are computed in \cite{TL}. Using that, we get the following table 
 	\begin{center}
 		\begin{tabular}{|c|c|c|c|c|} 
 			\hline
 			$(d,r,s, \delta)$ & $(3,11,0,1)$ & $(4,14,0,3)$ & $(5,17,0,6)$ & $(6,20,0,10)$ \\ 
 			\hline
 			$N_d^{\mathrm{Node}, \delta}(r,s)$ & $12960$ & $4057340$ & $2487128120$ & $2681467886460$ \\ 
 			\hline 
 		\end{tabular}
 	\end{center}
 	Finally, let us denote the by $\textnormal{Red}^{\textnormal{Node}, \delta}_d(r,s)$ to be the number of reducible 
 	planar degree $d$ curves in $\mathbb{P}^3$ with $\delta$ (unordered) nodes 
 	intersecting $r$ generic lines and $s$ generic points. 
 	This number can be computed using \cite[Proposition $8.4$]{TL}. Using that, we get 
 	\begin{center}
 		\begin{tabular}{|c|c|c|c|c|} 
 			\hline
 			$(d,r,s, \delta)$ & $(3,11,0, 1)$ & $(4,14,0,3)$ & $(5,17,0,6)$ & $(6,20,0,10)$ \\ 
 			\hline
 			$\textnormal{Red}_d^{\mathrm{Node}, \delta}(r,s)$ & $0$ & $329420$ & $507798840$ & $917948423100$ \\ 
 			\hline 
 		\end{tabular}
 	\end{center}
 	We now note that in all the cases we have tabulated, 
 	\[N_d^{\mathbb{P}^3, \mathrm{Planar}}(r,s) =  N_d^{\mathrm{Node}, \frac{(d-1)(d-2)}{2}}(r,s) - \textnormal{Red}_d^{\mathrm{Node}, \frac{(d-1)(d-2)}{2}}(r,s).\]
 	This is positive evidence for the fact that T.~Laarakker's formula for $N^{\textnormal{Node}, \delta}_d(r,s)$ 
 	is actually enumerative when $d \geq 1 + [\frac{\delta}{2}]$ (as opposed to $d \geq \delta$, which is proved in \cite{TL}). 
 	We also note that when $d=7$ and $\delta = \frac{(d-1)(d-2)}{2} = 15$, the formula for $N_d^{\mathrm{Node}, \frac{(d-1)(d-2)}{2}}(r,s)$ 
 	is not expected to be enumerative because of an obvious geometric reasons. To see why, suppose $s=0$ and $r = 23$.   
 	Through the required $r=23$ lines, we can place a double line 
 	through $4$ lines and through the remaining $19$ lines we can place a quintic that intersects the line so that the double line and the quintic 
 	lie in a plane. This is a degenerate configuration, and hence $N_7^{\mathrm{Node}, 15}(23, 0)$ is not expected to be enumerative.\\ 
 	\hf \hf This is analogous to the case of counting $\delta$-nodal degree $d$ curves in $\mathbb{P}^2$; let $N^{\delta}_d$ denote that number. 
 	A formula for this number can be  
 	explicitly found in \cite{FB} for instance. On the 
 	other hand, let $N_d^{\mathbb{P}^2}$ denote the number of rational degree $d$ curves in $\mathbb{P}^2$ through $3d-1$ generic points. 
 	Till $d=6$, we can verify that $N_d^{\mathbb{P}^2}$ is logically consistent with the corresponding value of $N^{\frac{(d-1)(d-2)}{2}}_d$ (after subtracting 
 	the number of irreducible curves). From $d=7$, we can not make any such consistency check, because $N^{15}_7$ is not enumerative; there are double 
 	lines that can pass through two of the $20$ points and a quintic through the remaining $18$ points. We also note that this fact is consistent with the  
 	{G}\"ottsche threshold of when the number $N^{\delta}_{d}$ is supposed to be enumerative; this is proved in \cite{KlSh}. 
 	Our computations give evidence 
 	to show that a similar bound is likely to be true for the case of planar curves in $\mathbb{P}^3$. \\
 	
 	
 	\section{Acknowledgement} 
 	The first author is grateful to Martijn Kool and Ties Laarakker for a very fruitful exchange of ideas and discussions 
 	related to the parallel question of counting $\delta$-nodal planar curves in $\mathbb{P}^3$, that is investigated in 
 	\cite{TL}. We are also very grateful  to ICTS for their hospitality and conducive atmosphere for doing mathematics research. 
 	We would like to specially acknowledge two programs conducted by ICTS where a significant part of the project was carried  out:  
 	Analytic and Algebraic Geometry (Code: ICTS/Prog-AAG2018/03) and 
 	Integrable Systems in Mathematics, Condensed Matter and Statistical Physics (Code: ICTS/integrability2018/07). The first author 
 	would also like to thank Center for Quantum Geometry of Moduli Space at Aarhus, Denmark (DNRF95) for giving him a chance to spend $6$ weeks 
 	there, when the author  got the initial idea for this project; the visit was 
 	was mainly paid by the grant ``EU-IRSES Fellowship within FP7/2007-2013 under grant 
 	agreement number 612534, project MODULI - Indo European Collaboration on Moduli Spaces." 
 	Finally, the first author would like to acknowledge the External Grant 
 	he has obtained, namely 
 	MATRICS (File number: MTR/2017/000439) that has been sanctioned by the Science and Research Board (SERB).

 \end{document}